\documentclass[aos,citesort,dvips]{arximspdf}
\usepackage{dcolumn}
\usepackage{graphicx}

\doi{10.1214/10-AOS821}
\volume{38}
\issue{6}
\pubyear{2010}
\firstpage{3605}
\lastpage{3629}

\makeatletter

\newcolumntype{d}[1]{D{.}{.}{#1}}

\newproclaim{Assumption}{Assumption}
\newtheorem{theorem}{Theorem}
\newtheorem{Proposition}{Proposition}
\newtheorem{lem}{Lemma}

\newcommand{\argmax}{\arg\max}

\newcommand{\bLambda}{\bolds{\Lambda}}
\newcommand{\bSigma}{\bolds{\Sigma}}

\newcommand{\bC}{\mathbf{C}}
\newcommand{\bD}{\mathbf{D}}
\newcommand{\bE}{\mathbf{E}}
\newcommand{\bF}{\mathbf{F}}
\newcommand{\bG}{\mathbf{G}}
\newcommand{\bH}{\mathbf{H}}
\newcommand{\bI}{\mathbf{I}}
\newcommand{\bM}{\mathbf{M}}
\newcommand{\bP}{\mathbf{P}}
\newcommand{\bS}{\mathbf{S}}
\newcommand{\bV}{\mathbf{V}}
\newcommand{\bU}{\mathbf{U}}
\newcommand{\bX}{\mathbf{X}}
\newcommand{\bY}{\mathbf{Y}}
\newcommand{\bZ}{\mathbf{Z}}

\newcommand{\ba}{\mathbf{a}}
\newcommand{\br}{\mathbf{r}}

\newcommand{\bmu}{\bolds{\mu}}

\makeatother

\begin{document}
\begin{frontmatter}

\title{Convergence and prediction of principal component scores in high-dimensional settings}
\runtitle{Convergence and prediction of PC scores}

\begin{aug}
\author[A]{\fnms{Seunggeun} \snm{Lee}\corref{}\thanksref{t1}\ead[label=e1]{slee@bios.unc.edu}},
\author[A]{\fnms{Fei} \snm{Zou}\thanksref{t1,t2}\ead[label=e2]{fzou@bios.unc.edu}} and
\author[A]{\fnms{Fred A.} \snm{Wright}\thanksref{t1,t2}\ead[label=e3]{fwright@bios.unc.edu}}
\runauthor{S. Lee, F. Zou and F. A. Wright}
\affiliation{University of North Carolina at Chapel Hill}
\address[A]{Department of Biostatistics \\
University of North Carolina at Chapel Hill \\
3101 McGavran-Greenberg, CB 7420\\
Chapel Hill, North Carolina 27599\\
USA\\
\printead{e1}\\
\phantom{E-mail: }\printead*{e2}\\
\phantom{E-mail: }\printead*{e3}} 
\end{aug}

\thankstext{t1}{Supported in part by NIH Grant GM074175.}
\thankstext{t2}{Supported in part by the Carolina Environmental
Bioinformatics Center (EPA RD832720) and a Gillings Innovation Award.}

\received{\smonth{12} \syear{2009}}
\revised{\smonth{2} \syear{2010}}

%
\begin{abstract}
A number of settings arise in which it is of interest to predict
Principal Component (PC) scores for new observations using data from
an initial sample. In this paper, we demonstrate that naive approaches
to PC score prediction can be substantially biased toward 0 in the
analysis of large matrices. This phenomenon is largely related to known
inconsistency results for sample eigenvalues and eigenvectors as both
dimensions of the matrix increase. For the spiked eigenvalue model for
random matrices, we expand the generality of these results, and
propose bias-adjusted PC score prediction. In addition, we compute the
asymptotic correlation coefficient between PC scores from sample and
population eigenvectors. Simulation and real data examples from the
genetics literature show the improved bias and numerical properties of
our estimators.
\end{abstract}

%

\end{frontmatter}

\section{Introduction}\label{sec1}
Principal component analysis (PCA)~\cite{jolliffe2002pca} is one of the
leading statistical tools for analyzing multivariate data. It is
especially popular in genetics/genomics, medical imaging and
chemometrics studies where high-dimensional data is common. PCA is
typically used as a dimension reduction tool. A small number of top
ranked principal component (PC) scores are computed by projecting data
onto spaces spanned by the eigenvectors of sample covariance matrix,
and are used to summarize data characteristics that contribute most to
data variation. These PC scores can be subsequently used for data
exploration and/or model predictions. For example, in genome-wide
association studies (GWAS), PC scores are used to estimate ancestries
of study subjects and as covariates to adjust for population
stratification~\cite{price2006pca,patterson2006psa}. In gene expression
microarray studies, PC scores are used as synthetic ``eigen-genes'' or\vadjust{\goodbreak}
``meta-genes'' intended to represent and discover gene expression
patterns that might not be discernible from single-gene analysis
\cite{wall2003svd}.

Although PCA is widely applied in a number of settings, much of our
theoretical understanding rests on a relatively small body of
literature. Girshick~\cite{girshick1936principal} introduced the idea
that the eigenvectors of sample covariance matrix are maximum
likelihood estimators.
Here, a key concept in a population view of PCA is that the data arise
as $p$-variate values from a distinct set of $n$ independent samples.
Later, the asymptotic distribution of eigenvalues and eigenvectors of
the sample covariance matrix (i.e., the sample eigenvalues and
eigenvectors) were derived for the situation where $n$ goes to
infinity and $p$ is fixed~\cite{girshick1939sampling,Anderson1}.
With the development of modern high-throughput technologies, it is not
uncommon to have data where $p$ is comparable in size to $n$, or
substantially larger. Under the assumption that $p$ and $n$ grow at
the same rate, that is $p/n \rightarrow\gamma> 0$, there has been
considerable effort to establish convergence results for sample
eigenvalues and eigenvectors (see review~\cite{bai1999msa}). The
convergence of the sample eigenvalues and eigenvectors under the
``spiked population'' model proposed by Johnstone~\cite{johnstone2001}
has also been established~\cite{Baik2006,paul2007,nadler2008fsa}. For
this model, it is well known that the sample eigenvectors are not
consistent estimators of the eigenvectors of population covariance
(i.e., the population eigenvectors)
\cite{johnstone2007spc,paul2007,nadler2008fsa}. Furthermore, Paul~\cite{paul2007} has derived the degree of discrepancy in terms of the
angle between the sample and population eigenvectors, under Gaussian
assumptions for $0 < \gamma< 1$. More recently, Nadler
\cite{nadler2008fsa} has extended the same result to the more general
$\gamma
> 0 $ using a matrix perturbation approach.

These results 
have considerable potential practical utility in understanding the
behavior of PC analysis and prediction in modern datasets, for which
$p$ may be large. The practical goals of this paper focus primarily on
the prediction of PC scores for samples which were not included in the
original PC analysis. For example, gene expression data of new breast
cancer patients may be collected, and we might want to estimate their
PC scores in order to classify their cancer sub-type. The
recalculation of PCs using both new and old data might not be
practical. For example, if the application of PCs from gene expression
is used as a diagnostic tool in clinical applications. For GWAS
analysis, it is known that PC analysis which includes related
individuals tends to generate spurious PC scores which do not reflect
the true underlying population substructures. To overcome this problem,
it is common practice to include only one individual per family/sibship
in the initial PC analysis. Another example arises in cross-validation
for PC regression, in which PC scores for the test set might be derived
using PCA performed on the training set~\cite{jackson2005user}. For
all of these applications, the predicted PC scores for a new sample are
usually estimated in the ``naive'' fashion, in which the data vector of
the new sample is multiplied by the sample eigenvectors from the
original PC analysis.\vadjust{\goodbreak}
Indeed, there appears to
be relatively little recognition in the genetics or data mining
literature that this
approach may lead to misleading conclusions.

%
\begin{figure}[b]
\vspace*{-3pt}
\includegraphics{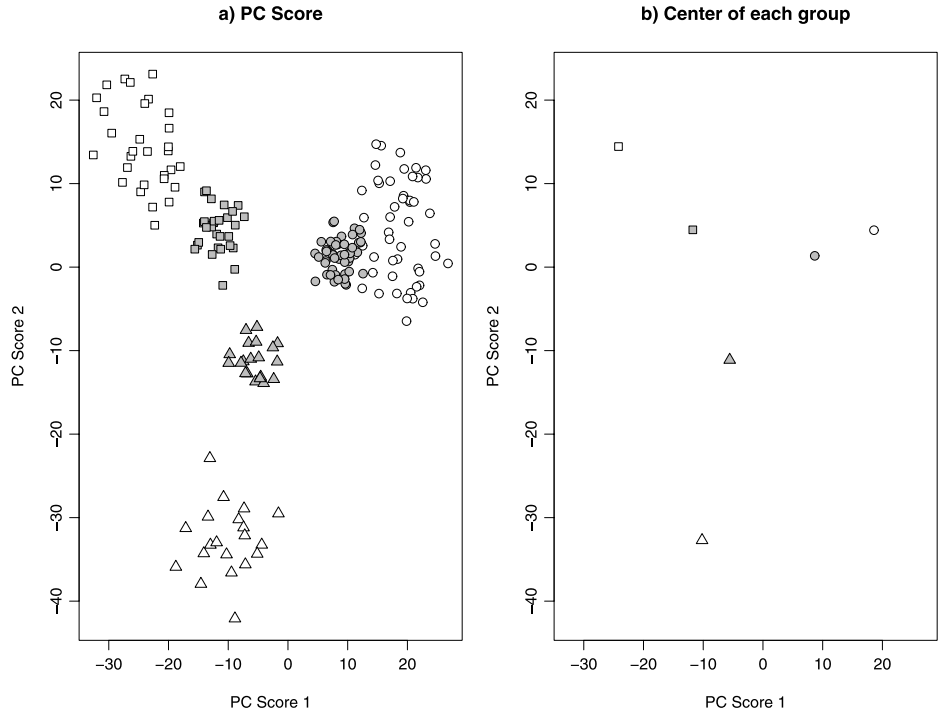}

\caption{Simulation results for $p=5000$ and $n=(50,30,20)$.
Different symbols represent different groups. White background color
represents the training set and grey background color represents the
test set. \textup{(a)} First 2 PC score plot of all simulated samples. \textup{(b)} Center
of each group.}\label{fig:1}
\end{figure}

For low-dimensional data, where $p$ is fixed as $n$ increases or
otherwise much smaller than $n$, the predicted PC scores are nearly
unbiased and well-behaved. However, for high-dimensional data,
particularly with $p>n$, they tend to be biased
and shrunken toward 0. The following simple example of a stratified
population with three strata
illustrates the shrinkage phenomenon for predicted PC scores. We
generated a training data set with
$n=100$ and $p=5000$. Among the 100 samples, 50 are from stratum 1, 30
are from stratum 2 and the rest from
stratum 3. For each stratum, we first created a $p$-dimensional mean
vector $\bmu_k$ $(k=1,2,3)$. Each element of each mean vector was
created by drawing randomly with replacement from $\{-0.3, 0, 0.3\}$,
and thereafter considered a fixed property of the stratum. Then for
each sample from the $k$th stratum, its $p$
covariates were simulated from the multivariate normal distribution
$\operatorname{MVN}(\bmu_k, 4 \bI)$, where $\bI$ is the $p
\times p$ identity matrix. A test dataset with the same sample size and
$\bmu_k$ vectors was also simulated.
Figure~\ref{fig:1}\vadjust{\goodbreak} shows that the predicted PC scores for the test data are much
closer to 0 compared to the scores from
the training data. This shrinkage phenomenon may create a serious
problem if the predicted PC scores are used to
classify new test samples, perhaps by similarity to previous apparent
clusters in the original data. In
addition, the predicted PC scores may produce incorrect results if used
for downstream analyses (e.g., as
covariates in association analyses).

In this paper, we investigate the degree of shrinkage bias associated
with the predicted PC scores, and then
propose new bias-adjusted PC score estimates. As the shrinkage
phenomenon is largely related to the limiting
behavior of the sample eigenvectors, our first step is to describe the
discrepancy between the sample and
population eigenvectors. To achieve this purpose, we follow the
assumption that $p$ and $n$ both are large and
grow at the same rate. By applying and extending results from random
matrix theory, we establish the
convergence of the sample eigenvalues and eigenvectors under the spiked
population model. We generalize Theorem
4 of Paul~\cite{paul2007}, which describes the asymptotic angle between
sample and population eigenvectors, to
non-Gaussian random variables for any $\gamma> 0$. We further derive
the asymptotic angle between PC scores
from sample eigenvectors and population eigenvectors, and the
asymptotic shrinkage factor of the PC score
predictions. Finally, we construct estimators of the angles and the
shrinkage factor.
The theoretical results are presented in Section~\ref{sec2}.

In Section~\ref{sec3}, we report simulations to assess the finite sample
accuracy of the proposed asymptotic angle and
shrinkage factor estimators. We also show the potential improvements in
prediction accuracy for PC regression by
using the bias-adjusted PC scores. In Section~\ref{sec4}, we apply our PC
analysis to a real genome-wide association
study, which demonstrates that the shrinkage phenomenon occurs in real
studies and that adjustment is needed.

\section{Method}\label{sec2}\vspace*{-2pt}
\subsection{General setting}
Throughout this paper, we use $^T$ to denote matrix transpose,
$\stackrel{p}{\rightarrow}$
to denote convergence in probability, and
$\stackrel{\mathrm{a.s.}}{\rightarrow}$ to denote almost sure convergence. Let
$\bLambda=
\operatorname{diag}(\lambda
_1,\lambda_2,\ldots,\lambda_p)$, a $p\times
p$ matrix with $\lambda_1\geq\lambda_2 \geq\cdots\geq\lambda_p$,
and $\bE=[\mathbf{e}_1,\ldots,\mathbf{e}_p]$, a $p\times p$
orthogonal matrix.

Define the $p\times n$ data matrix, $\bX$ as $[ \mathbf{x}_{1},\ldots
, \mathbf{x}
_{n}]$, where $\mathbf{x}_j$ is the
$p$-dimension\-al vector corresponding to the $j$th sample. For the
remainder of the paper, we assume the following.\vspace*{-2pt}
\begin{Assumption}\label{assum1}
$\bX= \bE\bLambda^{1/2}\bZ$, where $\bZ
=\{ z_{ij} \}$ is a $p \times n$ matrix
whose elements $z_{ij}$'s are i.i.d. random variables with
$ E(z_{ij}) = 0, E(z_{ij}^2) = 1$ and $E(z_{ij}^4) < \infty$.\vspace*{-2pt}
\end{Assumption}

Although the $z_{ij}$'s are i.i.d., Assumption~\ref{assum1} allows for very
flexible covariance structures for $\mathbf X$, and thus the results of
this paper are quite general.
The population\vadjust{\goodbreak} covariance matrix of $\bX$ is $\bSigma= \bE\bLambda
\bE^T$.
The sample covariance matrix $\bS$ equals
\[
\bS= \bX\bX^T / n = \bE\bLambda^{1/2} \bZ\bZ^T
\bLambda^{1/2} \bE^T /n.
\]

The $\lambda_k$'s are the underlying population eigenvalues. The spiked
population model defined in
\cite{johnstone2001} assumes that all the population eigenvalues are 1,
except the first $m$ eigenvalues. That
is, $\lambda_1 \geq\lambda_2 \geq\cdots\geq\lambda_m > \lambda
_{m+1}=\cdots= \lambda_p =1$.
The spectral decomposition of the sample covariance matrix is
\[
\bS= \bU\bD\bU^T ,
\]
where $\bD= \operatorname{diag}(d_1 ,d_2 ,\ldots,d_p) $ is a diagonal matrix of the
ordered sample eigenvalues and $\bU=
[\mathbf{u}_1,\ldots,\mathbf{u}_p]$ is the corresponding $p \times
p$ sample
eigenvector matrix. Then the PC score matrix is $\bP=[\mathbf{p}_1,
\mathbf{p}_2,
\ldots, \mathbf{p}_n ]$,
where $\mathbf{p}_{v}^T=\mathbf{u}_{v}^T \bX$ is the $v$th sample PC
score. For
a new observation $\mathbf{x}_{\mathrm{new}}$, its predicted PC score is similarly
defined as $\bU^T \mathbf{x}_{\mathrm{new}}$ with the $v$th (PC) score equal to
$q_v=\mathbf{u}_v^{T}\mathbf{x}_{\mathrm{new}}$.

\subsection{Sample eigenvalues and eigenvectors}

Under the classical setting of fixed~$p$, it is well known that the
sample eigenvalues and eigenvectors are
consistent estimators of the corresponding population eigenvalues and
eigenvectors~\cite{anderson}. Under the
``large $p$, large $n$'' framework, however, 
the consistency is not guaranteed.
The following two lemmas summarize and extend some known
convergence results.
%
\begin{lem}\label{lemma1}
Let $ p/n \rightarrow\gamma\geq0 $ as $n \rightarrow\infty$.

\begin{longlist}
\item When $\gamma=0$,
%
%
\begin{equation}
d_v \stackrel{\mathit{a.s.}}{\rightarrow}
\cases{\lambda_v, &\quad for $v \leq m$, \cr
1, &\quad for $v > m$.}
\end{equation}

\item  When $\gamma> 0$,
%
%
\begin{equation}\label{baik_1}
d_v \stackrel{\mathit{a.s.}}{\rightarrow}
\cases{\rho(\lambda_v), &\quad for $v \leq k$,\cr
\bigl(1+\sqrt{\gamma}\bigr)^2, &\quad for $v = k + 1$,}
\end{equation}
where $k$ is the number of $\lambda_v$ greater than $1+\sqrt{\gamma}$,
and $\rho(x) = x(1+\gamma/(x-1)) $.
\end{longlist}
\end{lem}

The result in (ii) is due to Baik and Silverstein~\cite{Baik2006}, while
the proof of (i) can be found in Section~\ref{L1}.
The result in (i) shows that when $\gamma=0$, the sample eigenvalues
converge to the corresponding population
eigenvalues, which is consistent with the classical PC result where $p$
is fixed. The result in (ii) shows that
for any nonzero $\gamma$, $d_v$ is no longer a consistent estimator of
$\lambda_v$. However, a consistent
estimator of $\lambda_v$ can be constructed from (\ref{baik_1}). Define
\[
\rho^{-1}(d) = \frac{ d + 1 -\gamma+ \sqrt{ (d + 1 -\gamma)^2 - 4 d}}{2}.
\]
Then $\rho^{-1}(d_v)$ is a consistent estimator of $\lambda_v$ when $
\lambda_v > 1+ \sqrt{\gamma}$.
Furthermore, Baik, Ben Arous and P{\'e}ch{\'e}~\cite{baik2005ptl} have
shown the $\sqrt{n}$-consistency of $d_v$ to $\rho(\lambda_v)$, and
Bai and Yao~\cite{bai2008clt} have shown that $d_v$ is asymptotically
normal.\vspace*{-2pt}
\begin{lem}\label{lemma2}
Suppose $ p/n \rightarrow\gamma\geq0 $ as $n \rightarrow\infty$.
Let $\langle\cdot,\cdot\rangle $ be an inner product between
two vectors. Under the
assumption of multiplicity one:

\begin{longlist}
\item if $ 0<\gamma<1 $, and the $z_{ij}$'s follow the standard normal
distribution, then
%
%
\begin{equation}\label{paul_1}
| \langle \mathbf{e}_v,\mathbf{u}_v \rangle | \stackrel{\mathit{a.s.}}{\rightarrow}
\cases{\phi(\lambda_v), &\quad if $\lambda_v > 1 + \sqrt{\gamma}$,\cr
0, &\quad if $1 < \lambda_v \leq1+\sqrt{\gamma}$;}
\end{equation}

\item removing the normal assumption on the $z_{ij}$'s, the following
weaker convergence result
holds for all $\gamma\geq0$:
%
%
\begin{equation}\label{eigen_angle}
| \langle \mathbf{e}_v,\mathbf{u}_v\rangle | \stackrel{p}{\rightarrow}
\cases{\phi(\lambda_v), &\quad if $\lambda_v > 1+\sqrt{\gamma}$, \cr
0, &\quad if $1 < \lambda_v \leq1+\sqrt{\gamma}$.}
\end{equation}
Here $ \phi(x) = \sqrt{(1-\frac{\gamma}{(x-1)^2}) / (1+\frac
{\gamma
}{x-1})} $.\vspace*{-2pt}
\end{longlist}
\end{lem}

The inner product between unit vectors is the cosine angle between
these two. Thus, Lemma~\ref{lemma2} shows the
convergence of the angle between population and sample eigenvectors.
For (i), Paul~\cite{paul2007} proved it for
$\gamma<1$; while Nadler~\cite{nadler2008fsa} obtained the same
conclusion for $\gamma> 0$ using the matrix
perturbation approach under the Gaussian random noise model. We relax
the Gaussian assumption on $z$ and prove
(ii) for $\gamma\geq0$ in Section~\ref{L2}. The result of (ii) is
general enough for the application of PCA to,
for example, genome-wide association mapping, where each entry of $\bX$
is a standardized variable of SNP
genotypes, which are typically coded as $\{0, 1, 2\}$, corresponding to
discrete genotypes.\vspace*{-2pt}

\subsection{Sample and predicted PC scores}

In this section, we first discuss convergence of the sample PC scores,
which forms the basis for the
investigation of the shrinkage phenomenon of the predicted PC scores.
For the sample PC scores, we have
the following theorem.\vspace*{-2pt}
\begin{theorem}\label{theo1}
Let $ \mathbf{g}_v^T = \mathbf{e}_v^T \bX/
\sqrt{\mathbf{n}\lambda_v} $, the normalized $v$th PC score derived
from a corresponding population eigenvector,
$\mathbf{e}_v$, and $\tilde{\mathbf{p}}_v = \mathbf{p}_v / \sqrt
{\mathbf{n}d_v}$, the
normalized $v$th sample PC score. Suppose $ p/n
\rightarrow\gamma\geq0 $ as $n \rightarrow\infty$. Under the
multiplicity one assumption,
%
%
\begin{equation}\label{score_angle}
| \langle \mathbf{g}_v,\tilde{\mathbf{p}}_v\rangle | \stackrel{p}{\rightarrow}
\cases{\sqrt{1-\dfrac{\gamma}{(\lambda_v-1)^2}}, &\quad if $\lambda_v >
1+\sqrt{\gamma}$,\vspace*{2pt}\cr
0, &\quad if $1 < \lambda_v \leq1+\sqrt{\gamma}$.}\vspace*{-2pt}
\end{equation}
\end{theorem}

The proof can be found in Section~\ref{T1}. In PC analysis, the sample
PC scores are typically used to estimate certain latent variables
(largely the PC scores\vadjust{\goodbreak} from population eigenvectors) that represent the
underlying data structures. The above result allows us to quantify the
accuracy of the sample PC scores. Note that here $\langle
\mathbf{g}_v,\tilde{\mathbf{p}}_v\rangle $ is the correlation coefficient between
$\mathbf{g}_v$ and $\tilde{\mathbf{p}}_v$. Compared to (\ref{paul_1})
in Lemma~\ref{lemma2}, the angle between the PC scores is smaller than
the angle
between their corresponding eigenvectors.

Before we formally derive the asymptotic shrinkage factor for the
predicted PC scores, we first describe in mathematical terms
the shrinkage phenomenon that was demonstrated in the \hyperref
[sec1]{Introduction}.
Note that the first population eigenvector $\mathbf{e}_1$ satisfies
\[
\mathbf{e}_1 = \mathop{\argmax}_{\ba\dvtx\ba^{T}\ba= 1 } E ( (\ba^{T}
\mathbf{x})^2 )
\]
for a random vector $\mathbf{x}$ that follows the same distribution of
the $\mathbf{x}
_j$'s. For the data
matrix $\bX$, its first sample eigenvector $\mathbf{u}_1$ satisfies
\[
\mathbf{u}_1 = \mathop{\argmax}_{\ba\dvtx\ba^{T}\ba= 1 } \sum_{j=1}^{n}
(\ba^{T} \mathbf{x}_j)^2.
\]
Assuming that $\mathbf{u}_1$ and the new sample $\mathbf{x}_{\mathrm{new}}$
are independent of
each other, we have
%
%
\begin{eqnarray}\label{shrink_e3}
E( (\mathbf{u}_1^T \mathbf{x}_{\mathrm{new}})^2 )  &=& E( E (\mathbf{u}_1^T
\mathbf{x}_{\mathrm{new}}
\mathbf{x}_{\mathrm{new}}^T \mathbf{u}_1^T |\mathbf{u}_1) )
= E( \mathbf{u}_1^T E (\mathbf{x}_{\mathrm{new}} \mathbf{x}_{\mathrm{new}}^T ) \mathbf
{u}_1^T ) \nonumber\\[-8pt]\\[-8pt]
&=& E( \mathbf{u}_1^T \bSigma\mathbf{u}_1^T ) \leq\mathbf{e}_1^T
\bSigma\mathbf{e}_1 = E( (\mathbf{e}
_1^T \mathbf{x}_{\mathrm{new}})^2 ).\nonumber
\end{eqnarray}
Since the $\mathbf{u}_1^T \mathbf{x}_j$'s $(j=1,\ldots,n)$ follow
the same distribution,
%
%
\begin{equation}\label{shrink_e4}\quad
n E((\mathbf{e}_1^T \mathbf{x}_j )^2) = E\Biggl( \sum_{j=1}^{n}(\mathbf
{e}_1^T \mathbf{x}_j )^2 \Biggr)
\leq E\Biggl( \sum_{j=1}^{n}(\mathbf{u}_1^T \mathbf{x}_j )^2 \Biggr) = n
E((\mathbf{u}_1^T \mathbf{x}_j )^2 ).
\end{equation}
By (\ref{shrink_e3}) and (\ref{shrink_e4}), we can show that
\[
E( (\mathbf{u}_1^T \mathbf{x}_{\mathrm{new}})^2 ) \leq E( (\mathbf{e}_1^T
\mathbf{x}_{\mathrm{new}})^2 )
= E( (\mathbf{e}_1^T \mathbf{x}_j )^2 ) \leq E( (\mathbf{u}_1^T
\mathbf{x}_j )^2 ),
\]
%
which demonstrates the shrinkage feature of the predicted PC scores.
The amount of the shrinkage, or the
asymptotic shrinkage factor, is given by the following theorem.
\begin{theorem}\label{theo2}
Suppose $ p/n \rightarrow\gamma\geq0 $ as
$n \rightarrow\infty$, $\lambda_v > 1 + \sqrt{\gamma}$. Under the
multiplicity one assumption,
%
%
\begin{equation}\label{shrink}
\sqrt{\frac{ E( q_v^2)}{ E( p_{vj}^2) }}
\stackrel{n\rightarrow\infty}{\rightarrow}
\frac{\lambda_v -1}{\lambda_v +\gamma-1},
\end{equation}
where $p_{vj}$ is the $j$th element of $\mathbf{p}_{v}$.
\end{theorem}

The proof is given in Section~\ref{T2}. We call $ (\lambda_v
-1)/(\lambda_v +\gamma-1)$, the (asymptotic)
shrinkage factor for a new subject.
As shown, the shrinkage factor is smaller than $1$ if $\gamma> 0$.
Quite sensibly, it is a decreasing function of
$\gamma$ and an increasing function of $\lambda_v$. The bias of the
predicted PC score can be
potentially large\vadjust{\goodbreak} for those high-dimensional data where $p$ is
substantially greater than $n$, and/or for the
data with relatively minor underlying structures where $\lambda_v$ is small.

\subsection{Rescaling of sample eigenvalues}

The previous theorems are based on the assumption that all except the
top $m$ eigenvalues are equal to 1.
Even under the spiked eigenvalue model, some rescaling of the sample
eigenvalues may be necessary with real data.

For a given data, let its ordered population eigenvalues $ \bLambda^{*}
= \{ \zeta\lambda_1 ,\ldots,\zeta\lambda_m$, $\zeta,\ldots,\zeta\}$,
where $\zeta\neq1$, and its corresponding sample eigenvalues $ \bD
^{*} = \{
d_1^{*} ,\ldots,\break d_n^{*} \} $. We can show that (\ref{eigen_angle}),
(\ref{shrink}) and
(\ref{score_angle}) still hold under such circumstances. However,
$\rho
^{-1}(d_v^{*})$ is no longer a consistent
estimator of $\lambda_v$, because
\[
d_v^{*} \stackrel{\mathrm{a.s.}}{\rightarrow} \zeta\lambda_v\biggl(1+\frac{\gamma
}{\lambda_v-1}\biggr)= \zeta\rho(\lambda_v).
\]
To address this issue, Baik and Silverstein~\cite{Baik2006} have
proposed a simple approach to estimate $\zeta$. In their method,
the top significant large sample eigenvalues are first separated from
the other grouped sample eigenvalues. Then
$\zeta$ is estimated as the ratio between the average of the grouped
sample eigenvalues and the mean determined by the
Mar{\v{c}}enko--Pastur law~\cite{marvcenko1967}. To separate the eigenvalues,
they have suggested to use a screeplot of the percent variance versus
component number. However, for
real data, we may not be able to clearly separate the sample
eigenvalues in such a manner and readily apply the approach. Thus, we
need an automated
method which does not require a clear separation of the sample eigenvalues.

The expectation of the sum of the sample eigenvalues when $\zeta= 1$ is
\[
E\Biggl(\sum_{v=1}^p d_v\Biggr) = E(\operatorname{trace}(\bS)) = \operatorname{trace}(E(\bS))
= \operatorname{trace}(\bSigma) =
\sum_{v=1}^p \lambda_v .
\]
Thus, the sum of the rescaled eigenvalues is expected to be close to
$(\sum_{v=1}^m \lambda_v + p - m)$.
Let $ r_v = d_v^{*}/ ( \sum_{v=1}^{p} d_v^{*} ) $ and $ \hat{d}_{v} $
be a properly rescaled
eigenvalue, then $ \hat{d}_{v} $ should be very close to $ r_v (\sum
_{v=1}^m \lambda_v + p - m) $. Note
that ${p}/(\sum_{v=1}^m \lambda_v + p -\break m) \rightarrow1$ for fixed
$m$ and $\lambda_v$. Thus, $ p r_v $
is a properly adjusted eigenvalue. However, for finite $n$ and $p$, the
difference between $p$ and
$(\sum_{v=1}^m \lambda_v +p -m) $ can be substantial, especially when
the first several $\lambda_v$'s are
considerably larger than $1$. To reduce this difference, we propose a
novel method which
iteratively estimates the $(\sum_{v=1}^m \lambda_v + p - m)$ and
$\hat{d}_{v}$.



1. Initially set $\hat{d}_{v,0} = p r_v $.

2. For the $l$th iteration, set $\hat{\lambda}_{v,l} = \rho
^{-1}(\hat
{d}_{v,l-1}) $ for
$ \hat{d}_{v,l-1} > (1+\sqrt{\gamma})^2 $, and $\hat{\lambda
}_{v,l} =
1$ for $ \hat{d}_{v,l-1}
\leq(1+\sqrt{\gamma})^2 $. Define
$k_l$ as the number of $\hat{\lambda}_{v,l}$'s that are greater than 1,
and let
\[
\hat{d}_{v,l} = \Biggl( \sum_{v=1}^{k_l} \hat{\lambda}_{v,l} + p - k_l \Biggr)
r_v .\vadjust{\goodbreak}
\]

3. If $\sum_{v=1}^{k_l} \hat{\lambda}_{v,l} + p - k_l$ converges, let
\[
\hat{d}_v = \hat{d}_{v,l}
\]
and stop.
Otherwise, go to step 2.

The consistency of $\hat{d}_v$ to $\rho(\lambda_v)$ is shown in the
following theorem.
\begin{theorem}\label{theo3}
Let $\hat{d}_v$ be the rescaled sample eigenvalue from the proposed
algorithm. Then, for $\lambda_v
> 1 + \sqrt{\gamma}$ with multiplicity one,
\[
\hat{d}_v \stackrel{p}{\rightarrow} \rho(\lambda_v).
\]
\end{theorem}

Since $ \rho^{-1}(\hat{d}_v) \stackrel{p}{\rightarrow}\lambda_v $,
$\phi(\rho
^{-1}(\hat
{d}_v))^2$ is a
consistent estimator of $\phi(\lambda_v)^2$. Combining this fact with
Theorems~\ref{theo1} and~\ref{theo2}, we can obtain the
bias-adjusted PC score $q_v^{*}$
\[
q^{*}_v = q_v\frac{\rho^{-1}(\hat{d}_v) + \gamma-1}{\rho^{-1}(\hat
{d}_v) -1}
\]
and the asymptotic correlation coefficient between $\mathbf{g}_v$ and
$\mathbf
{\tilde{p}}_v $
\[
\sqrt{\biggl(1-\frac{\gamma}{(\rho^{-1}(\hat{d}_v)-1)^2}\biggr)}.
\]

\section{Simulation}\label{sec3}

\begin{figure}

\includegraphics{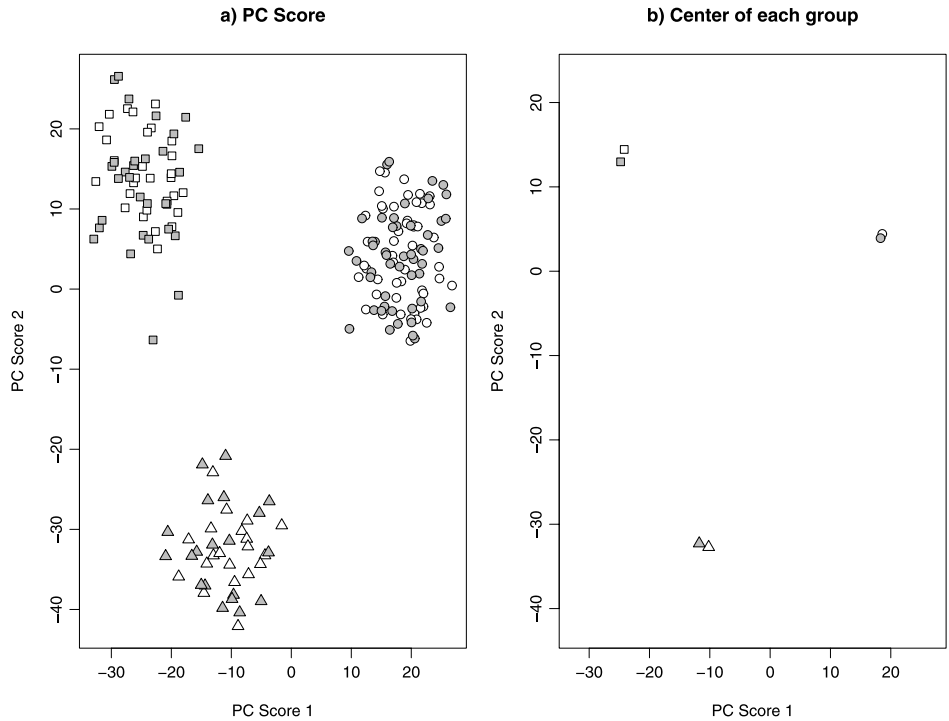}

\caption{Shrinkage-adjusted PC scores of the data in Figure
\protect\ref{fig:1}. Different symbols represent different groups. White background color
represents the training set and grey background color represents the
test set. \textup{(a)} Plots of all simulation samples. \textup{(b)} Center of each
group.}\label{fig:2}\vspace*{3pt}
\end{figure}

First, we applied our bias-adjustment process to the simulated data
described in the \hyperref[sec1]{Introduction}. Our
estimated asymptotic shrinkage factors are 0.465 and 0.329 for the
first and second PC scores, respectively. The
scatter plot of the top two bias-adjusted PC scores is given in Figure
\ref{fig:2}. After the bias adjustment,
%
%
the predicted PC scores of the test data are comparable to those of the
training data. This indicates that our
method is effective in correcting for the shrinkage bias.

Next, we conducted a new simulation to check the accuracy of our
estimators. For the $j$th sample ($j=1,\ldots,
n$), its $i$th variable was generated~as
\[
x_{ij} = \cases{\lambda_1 z_{ij}, &\quad $i = 1$, \cr
\lambda_2 z_{ij}, &\quad $i = 2$, \cr
z_{ij}, &\quad $i > 2$,}
\]
where $ \lambda_1 > \lambda_2 > 1$ and $z_{ij} \sim N(0,2^2)$. Under
this setting, $\lambda_1$ and $\lambda_2$
are the first and the second population eigenvalues. The first and
second population eigenvectors are $e_1 =
\{1,0,\ldots,0 \}$ and $e_2 = \{0,1,0,\ldots,0 \}$, respectively. We
set the standard deviation of $z_{ij}$ to 2 instead of
1, which allows us to test whether the rescaling procedure works
properly. We tried different values of
$\gamma$ and $n$, but set $\lambda_1$ and $\lambda_2$ to $4(1+\sqrt
{\gamma})$ and $2(1+\sqrt{\gamma})$,
respectively.

We split the simulated samples into test and training sets, each with
$n$ samples. We first estimated
the asymptotic shrinkage factor based on the training samples. We then
calculated the predicted PC scores on\vadjust{\goodbreak}
the test samples. To assess the accuracy of shrinkage factor estimator
for each PC, we empirically estimated the
shrinkage factor by the ratio of the mean predicted PC scores of the
test samples to the mean PC scores of the
training samples. That is, for the $v$th
PC, the
empirical shrinkage factor is estimated by $\sqrt{ \sum
_{i=1}^{n}q_{vi}^2 /\sum_{k=1}^{n}p_{vk}^2
}$.
On the training samples, we also estimated the empirical angle between
the sample and (known) population
eigenvectors, as well as the empirical angle between PC scores from
sample and population eigenvectors. The
asymptotic theoretical estimates were also calculated.
Tables~\ref{table:1} and~\ref{table:2} summarize the simulation results.
Our asymptotic estimators provide accurate estimates for the angles
and the shrinkage factor.\looseness=1

%
%
\begin{table}
\tabcolsep=0pt
\caption{Cosine angle estimates of eigenvectors and PC scores based on 1000
simulations. ``Angle'' indicates the~theoretical asymptotic
$\cos(\mathrm{angle})$, ``Estimate1'' indicates the empirical $\cos(\mathrm{angle})$
estimator, ``Estimate2'' indicates the asymptotic $\cos(\mathrm{angle})$ estimator.
For each estimator, each entry represents mean of $1000$
simulation results~with standard error in parentheses}
\label{table:1}
\begin{tabular*}{\tablewidth}{@{\extracolsep{\fill}}lccccccc@{}}
\hline
& & \multicolumn{3}{c}{\textbf{PC 1}} & \multicolumn{3}{c@{}}{\textbf{PC 2}}\\[-4pt]
& & \multicolumn{3}{c}{\hrulefill} & \multicolumn{3}{c@{}}{\hrulefill}\\
& & & \textbf{Angle} & \textbf{Angle} & & \textbf{Angle} & \textbf{Angle} \\
$\bolds\gamma$ &$\bolds n$& \textbf{Angle} & \textbf{Estimate1}
& \textbf{Estimate2} & \textbf{Angle} & \textbf{Estimate1} &
\textbf{Estimate2} \\
\hline
\multicolumn{8}{l}{Eigenvectors} \\[2pt]
\phantom{00}$ 1 $ &$ 100 $ & $0.93$ & $0.93\ (0.013)$ & $ 0.91\ (0.027) $ &$0.82$&
$ 0.81\ (0.053) $& $ 0.80\ (0.052) $ \\
$ $ &$ 200 $ && $0.93\ (0.009)$ & $ 0.92\ (0.014) $ && $ 0.81\ (0.030)
$& $ 0.81\ (0.032) $ \\[2pt]
\phantom{0}$ 20 $ & $ 100 $ & $ 0.70 $ & $ 0.69\ (0.037) $ & $ 0.70\ ( 0.031 ) $ &
$ 0.51 $& $ 0.50\ (0.053) $& $ 0.50\ (0.058) $ \\
& $ 200 $ & & $ 0.69\ (0.023) $ & $ 0.70\ ( 0.022 ) $ && $ 0.51\
(0.036) $& $ 0.51\ (0.041) $ \\[2pt]
$ 100 $ & $ 100 $ & $ 0.53 $ & $ 0.53\ (0.034) $ & $ 0.53\ ( 0.031 ) $
& $ 0.37 $& $ 0.35\ (0.043) $& $ 0.35\ (0.047) $ \\
& $ 200 $ && $ 0.53\ (0.024) $ & $ 0.53\ ( 0.024 ) $ && $ 0.36\
(0.029) $& $ 0.36\ (0.033) $ \\[2pt]
$ 500 $ & $ 100 $ & $ 0.38 $ & $ 0.38\ (0.029) $ & $ 0.38\ ( 0.028 ) $
& $ 0.25 $& $ 0.24\ (0.033) $& $ 0.24\ (0.037) $ \\
& $ 200 $ & & $ 0.38\ (0.020) $ & $ 0.38\ ( 0.020 ) $ && $ 0.25\
(0.021) $& $ 0.25\ (0.024) $ \\
[2pt]
\multicolumn{8}{l}{PC scores} \\
[2pt]
\phantom{00}$ 1 $ & $ 100 $ & $0.99$& $0.99\ (0.004)$ & $ 0.98\ (0.016) $ & $0.94$& $
0.93\ (0.036) $& $ 0.91\ (0.048) $ \\
$ $ & $ 200 $ && $0.99\ (0.003)$ & $ 0.99\ (0.006) $ && $ 0.94\ (0.019)
$& $ 0.93\ (0.024) $ \\[2pt]
\phantom{0}$ 20 $ & $ 100 $ & $ 0.98 $ & $ 0.97\ (0.083) $ & $ 0.98\ ( 0.008 ) $ &
$ 0.89 $& $ 0.86\ (0.105) $& $ 0.87\ (0.055) $ \\
& $ 200 $ & & $ 0.97\ (0.055) $ & $ 0.98\ ( 0.005 ) $ && $ 0.88\
(0.073) $& $ 0.88\ (0.036) $ \\[2pt]
$ 100 $ & $ 100 $ & $ 0.97 $ & $ 0.97\ (0.079) $ & $ 0.97\ ( 0.009 ) $
& $ 0.88 $& $ 0.85\ (0.109) $& $ 0.86\ (0.060) $ \\
& $ 200 $ && $ 0.97\ (0.058) $ & $ 0.97\ ( 0.006 ) $ && $ 0.86\
(0.076) $& $ 0.87\ (0.039) $ \\[2pt]
$ 500 $ & $ 100 $ & $ 0.97 $ & $ 0.96\ (0.084) $ & $ 0.97\ ( 0.010 ) $
& $ 0.87 $& $ 0.83\ (0.117) $& $ 0.84\ (0.069) $ \\
& $ 200 $ & & $ 0.96\ (0.058) $ & $ 0.97\ ( 0.007 ) $ && $ 0.86\
(0.076) $& $ 0.86\ (0.038) $ \\
\hline
\end{tabular*}
\end{table}

%
\begin{table}
\caption{Shrinkage factor estimates based on 1000 simulation. ``Factor''
indicates the theoretical asymptotic factor, ``Estimate1'' indicates
the empirical shrinkage factor estimator, ``Estimate2'' indicates the
asymptotic shrinkage factor estimator. For each estimator, each entry
represents mean of $1000$ simulation results with standard error in
parentheses}
\label{table:2}
\begin{tabular*}{\tablewidth}{@{\extracolsep{\fill}}lccccccc@{}}
\hline
& & \multicolumn{3}{c}{\textbf{PC 1}} & \multicolumn{3}{c@{}}{\textbf{PC 2}}\\[-4pt]
& & \multicolumn{3}{c}{\hrulefill} & \multicolumn{3}{c@{}}{\hrulefill}\\
& & & \textbf{Factor} & \textbf{Factor} && \textbf{Factor} & \textbf{Factor} \\
$\bolds\gamma$ & $\bolds n$ & \textbf{Factor} & \textbf{Estimate1} & \textbf{Estimate2}
& \textbf{Factor} & \textbf{Estimate1} & \textbf{Estimate2} \\
\hline
\phantom{00}$ 1 $ & $ 100 $ & 0.88& $0.88\ (0.017)$ & $ 0.87\ (0.076) $ &0.75& $
0.75\ (0.044) $& $ 0.76\ (0.063) $ \\
$ $ & $ 200 $ && $0.88\ (0.013)$ & $ 0.87\ (0.054) $ && $ 0.75\ (0.027)
$& $ 0.75\ (0.044) $ \\[2pt]
\phantom{0}$ 20 $ & $ 100 $ & $ 0.51 $ & $ 0.51\ (0.037) $ & $ 0.51\ ( 0.038 ) $ &
$ 0.33 $& $ 0.34\ (0.033) $& $ 0.32\ (0.038) $ \\
& $ 200 $ & & $ 0.51\ (0.025) $ & $ 0.51\ ( 0.026 ) $ && $ 0.34\
(0.022) $& $ 0.33\ (0.028) $ \\[2pt]
$ 100 $ & $ 100 $ & $ 0.30 $ & $ 0.30\ (0.024) $ & $ 0.30\ ( 0.030 ) $
& $ 0.17 $& $ 0.17\ (0.019) $& $ 0.17\ (0.023) $ \\
& $ 200 $ && $ 0.30\ (0.017) $ & $ 0.30\ ( 0.023 ) $ && $ 0.18\
(0.013) $& $ 0.17\ (0.017) $ \\[2pt]
$ 500 $ & $ 100 $ & $ 0.16 $ & $ 0.15\ (0.014) $ & $ 0.16\ ( 0.020 ) $
& $ 0.08 $& $ 0.08\ (0.010) $& $ 0.08\ (0.013) $ \\
& $ 200 $ && $ 0.15\ (0.010) $ & $ 0.16\ ( 0.014 ) $ && $ 0.08\
(0.007) $& $ 0.08\ (0.009) $ \\
\hline
\end{tabular*}
\end{table}

Finally, we conducted simulation to demonstrate an application of the
bias-adjusted PC scores in PC
regression. PC regression has been widely used in microarray
gene-expression studies~\cite{bovelstad2007}. In
this simulation, we let $p=5000$, and our set up is very similar to
the first simulation of
Bair et al.~\cite{bair2006prediction}. Let $x_{ij}$ denote the gene
expression level of the $i$th gene for the $j$th
subject. We generated each $x_{ij}$ according to
\[
x_{ij} = \cases{3 + \varepsilon, &\quad $i \leq g, j \leq n/2$, \cr
4 + \varepsilon, &\quad $i \leq g, j > n/2$, \cr
3.5 + \varepsilon, &\quad $i > g$,}\vadjust{\goodbreak}
\]
and the outcome variable $y_j$ as
\[
y_j = \frac{2}{g}\sum_{i=1}^g x_{ij} +\varepsilon_y,
\]
where $n$ is the number of samples, $g$ is the number of genes that are
differentially expressed and associated
with the phenotype, $\varepsilon\sim N(0,2^2)$ and $\varepsilon_y
\sim
N(0,1)$. A~total of eight different combinations of $n$ and $g$ were simulated.
For the training data,
we fit the PC regression with the first PC as the covariate and
computed the mean square error (MSE). For the
test samples with the same configuration of the training samples, we
applied the PC model built on the training
data to predict the phenotypes using the unadjusted and adjusted PC
scores. 
The results are presented in Table~\ref{table:3}. We see that the MSE of the test
%
%
\begin{table}[b]
\caption{Mean Square Error (MSE) of the PC regression based on gene-expression
microarray data simulation with and without shrinkage adjustment. 1000
simulation were conducted. Each entry in the table represents mean of
the MSE with standard error in parentheses}
\label{table:3}
\begin{tabular*}{\tablewidth}{@{\extracolsep{\fill}}ld{4.0}ccc@{}}
\hline
& & \textbf{Test data} & \textbf{Test data} & \\
$\bolds n$ & \multicolumn{1}{c}{$\bolds g$}
& \textbf{without adjustment} & \textbf{with adjustment}
& \textbf{Training data} \\
\hline
$ 100 $ & 150  & $1.97\ (0.256)$ & $ 1.70\ (0.284) $ & $ 1.61\
(0.284)$ \\
$ 100 $ & 300 & $1.63\ (0.230)$ & $ 1.17\ (0.167) $ & $ 1.12\
(0.158)$ \\
$ 100 $ & 500 & $1.43\ (0.204)$ & $ 1.07\ (0.157) $ & $ 1.03\
(0.147)$ \\
$ 100 $ & 1000 & $1.22\ (0.182)$ & $ 1.03\ (0.148) $ & $ 0.99\
(0.142)$ \\
$ 200 $ & 150 & $1.73\ (0.159)$ & $ 1.33\ (0.133) $ & $ 1.30\
(0.131)$ \\
$ 200 $ & 300 & $1.39\ (0.139)$ & $ 1.08\ (0.105) $ & $ 1.07\
(0.110)$ \\
$ 200 $ & 500 & $1.24\ (0.131)$ & $ 1.04\ (0.105) $ & $ 1.01\
(0.101)$ \\
$ 200 $ & 1000 & $1.10\ (0.114)$ & $ 1.02\ (0.101) $ & $ 1.00\
(0.101)$ \\
\hline
\end{tabular*}
\end{table}
set without bias adjustment is appreciably
higher than that of the test set with bias adjustment, and the MSE of
the test set with bias adjustment is
comparable with the MSE of the training~set.

\section{Real data example}\label{sec4}

Here, we demonstrate that the shrinkage phenomenon appears in real
data, and can be adjusted by our method. For
this purpose, genetic data on samples from unrelated individuals in the
Phase 3 HapMap study
(\url{http://hapmap.ncbi.nlm.nih.gov/}) were used. HapMap is a dense
genotyping study designed to elucidate
population genetic differences. The genetic data are discrete, assuming
the values 0, 1 or 2 at each genomic
marker (also known as SNPs) for each individual. Data from CEU
individuals (of northern and western European
ancestry) were compared with data from TSI individuals (Toscani
individuals from Italy, representing southern
European ancestry).

%
\begin{figure}[b]

\includegraphics{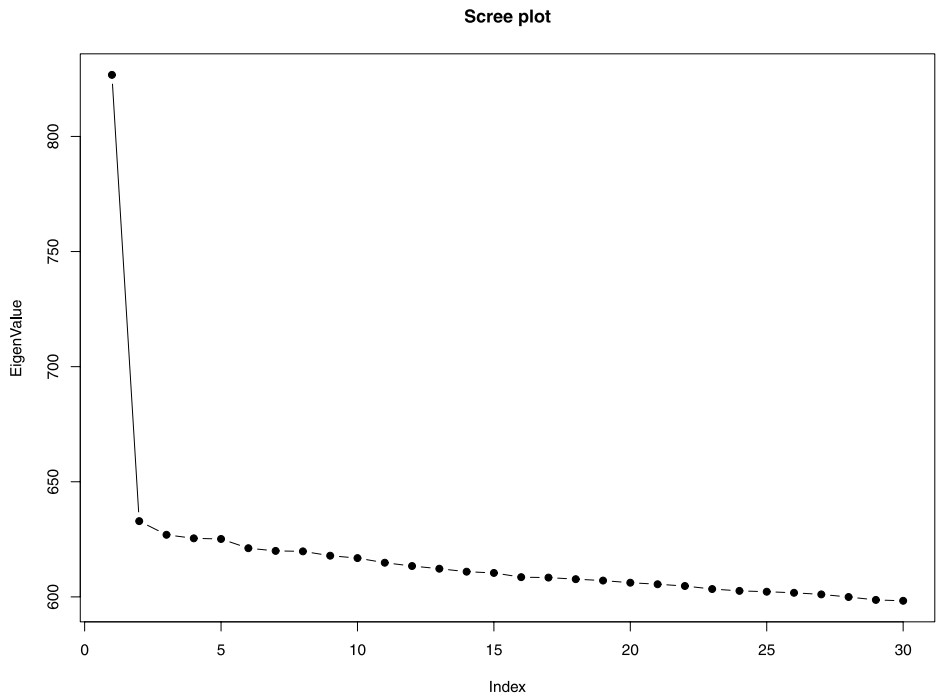}

\caption{Scree plot of the first 30 sample eigenvalues, CEU${}+{}$TSI
dataset.}\label{fig:3}
\end{figure}

Some initial data trimming steps are standard in genetic analysis.
We first removed apparently
related samples, and removed genomic markers with more than a 10\%
missing rate, and those with frequency less
than 0.01 for the minor genetic allele. To avoid spurious PC results,
we further pruned out SNPs that are
in high linkage disequlibrium (LD)~\cite{fellay2007}. Lastly, we
excluded $7$ samples with PC scores greater than
6 standard deviations away from the mean of at least one of the top
significant PCs [i.e., with Tracy--Widom (TW)
Test $p$-value $< $0.01]~\cite{price2006pca,patterson2006psa}. The final
dataset contained 178 samples (101 CEU, 77
TSI) and 100,183 markers. We mean-centered and variance-standardized
the genotypes for each marker
\cite{price2006pca}. The screeplot of the sample eigenvalues is
presented in Figure~\ref{fig:3}. The first
eigenvalue is substantially larger than the rest of the eigenvalues,
although the TW test actually identifies two
significant PCs. Figure~\ref{fig:3} suggests that our data approximately
satisfies the spiked eigenvalue assumption.

\begin{figure}

\includegraphics{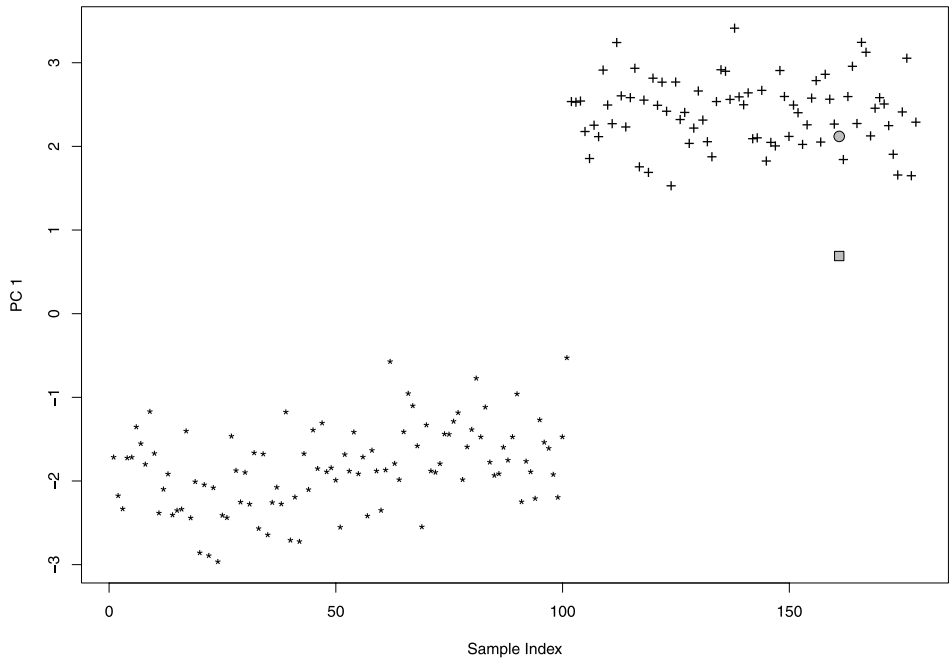}

\caption{An instance with and without shrinkage adjustment,
performed on Hapmap CEU(*) and TSI($+$). ``*'' and ``$+$'' represent PC
scores using all data. The $161$st sample was excluded from PCA,
and PC score for it was predicted. The grey rectangle represents the
predicted PC score without shrinkage adjustment and the grey circle
represents the predicted PC score after the shrinkage adjustment.}
\label{fig:4}
\end{figure}

We estimated the asymptotic shrinkage factor and compared it with the
following jackknife-based shrinkage factor estimate. For the first PC,
we first computed the scores of all samples. Next, we removed one
sample at a time and computed the (unadjusted) predicted PC score. We then
calculated the jackknife estimate as the square root
of the ratio of the means of the sample PC score and the predicted PC
score. The jackknife shrinkage factor
estimate is $0.319$, which is close to our asymptotic estimate $0.325$.
Figure~\ref{fig:4} shows the PC scores from the
%
%
whole sample, the predicted PC score of an illustrative excluded
sample, and its bias-adjusted predicted score.
Clearly, the predicted PC score without adjustment is very biased
toward zero, while the bias-adjusted PC
score is not.

\section{Discussion and conclusions}\label{sec5}

In this paper, we have identified and explored the shrinkage phenomenon
of the predicted PC scores, and have
developed a novel method to adjust these quantities. We also have
constructed the asymptotic estimator of
correlation coefficient between PC scores from population eigenvectors
and sample eigenvectors. In simulation
experiments and real data analysis, we have demonstrated the accuracy
of our estimates, and the capability to
increase prediction accuracy in PC regression by adopting shrinkage
bias adjustment. For achieving these, we
consider asymptotics in the large $p$, large $n$ framework, under the
spiked population model.

We believe that this asymptotic regime applies well to many
high-dimen\-sional datasets. It is not, however, the
only model paradigm applied to such data. For example, the large $p$
small $n$ paradigm
\cite{hall2005geometric,ahn2007high}, which assumes $p/n \rightarrow
\infty$, has also been explored. Under
this assumption, Jung and Marron~\cite{sungkyu} have shown that the
consistency and the strong inconsistency of the sample
eigenvectors to population eigenvectors depend on whether $p$ increases
at a slower or faster rate than
$\lambda_v$.
It may be argued that for real data where $p/n$ is ``large,'' we should
follow the paradigm of
Hall, Marron and Neeman~\cite{hall2005geometric}, Ahn et al. \cite
{ahn2007high}. However, for any real study, it is unclear
how to test whether $p$ increases at
a faster rate than $\lambda_v$, or vice versa, making the application
of Hall, Marron and Neeman~\cite{hall2005geometric}, Ahn et al. \cite
{ahn2007high} difficult in practice. Furthermore, the scenario where
$p$ and
$\lambda_v$ grow at the same rate is scientifically more interesting,
for which
we are aware of no theoretical results.
In contrast, our asymptotic results can be straightforwardly applied.
Further, our
simulation results indicate that for $p/n$ as large as 500, our
asymptotic results still hold well. We believe
that the approach we describe here applies to many datasets.

Although the results from the spiked model are useful, it is likely
that observed data has more structure than
allowed by the model. Recently, several methods have been suggested to
estimate population eigenvalues under more general scenarios \cite
{elkaroui2008sel,rao2008}. However, no analogous results are available
for the eigenvectors.
In data analysis, jackknife estimators, as demonstrated in the real
data analysis section, can be used.
However, resampling approaches are very computationally intensive, and
it remains of interest to establish the
asymptotic behavior of eigenvectors in a variety of situations.

We note that inconsistency of the sample eigenvectors does not
necessarily imply poor performance of PCA. For example, PCA has been
successfully applied in genome-wide association studies for accurate
estimation of ethnicity~\cite{price2006pca}, and in PC
regression for microarrays~\cite{ma2006additive}. However, for any
individual study we cannot rule out the possibility of poor performance
of the PC analysis. Our asymptotic result on the correlation
coefficient between PC scores from sample and population eigenvectors
provides us a measure to quantify the performance of PC analysis.

For the CEU/TSI data, SNP pruning was applied to adjust for strong LD
among adjacent SNPs. Such SNP pruning
is a common practice in the analysis of GWAS data, and has been
implemented in the popular GWAS analysis software
Plink~\cite{purcell2007pts}. The primary goal of SNP pruning is to
avoid spurious PC results\vadjust{\goodbreak}
unrelated to population substructures. Technically, our approach does
not rely on any independence assumption of the
SNPs. However, strong local correlation may affect eigenvalues
considerably. Thus, the value in SNP pruning may be viewed as helping
the data better accord with the assumptions of the spiked population model.
From the CEU/TSI data and our experience in other
GWAS data, we have found that the most common pruning procedure
implemented in Plink is sufficient for us to then apply our
methods.

\section{Proofs}\label{sec6}

Note that $\bE\bLambda^{1/2} \bZ\bZ^T \bLambda^{1/2} \bE^T$ and
$\bLambda^{1/2}\bZ\bZ^T \bLambda^{1/2}$ have
the same eigenvalues, and $\bE^T \bU$ is the eigenvector matrix of
$\bLambda^{1/2}\bZ\bZ^T \bLambda^{1/2}$.
Since eigenvalues and angles between sample and population eigenvectors
are what we concerned about, without
loss of generality (WLOG), in the sequel, we assume $\bLambda$ to be
the population covariance matrix.

\subsection{Notation}

We largely follow notation in Paul~\cite{paul2007}. We denote $\lambda
_v(\bS)$ as the $v$th largest eigenvalue
of $\bS$. Let suffice $A$ represent the first $m$ coordinates and $B$
represent the remaining coordinates. Then
we can partition $\bS$ into
\[
\bS= \left[
\matrix{
\bS_{AA} & \bS_{AB} \cr
\bS_{BA} & \bS_{BB} }
\right].
\]
We similarly partition the $v$th eigenvector $\mathbf{u}_v^T$ into
$(\mathbf{u}
_{A,v},\mathbf{u}_{B,v})$ and $\bZ^T$ into
$[\bZ_A^T , \bZ_B^T ] $. Define $ R_v$ as $\Vert\mathbf
{u}_{B,v}\Vert$ and let
$\ba_v = \mathbf{u}_{A,v}/\sqrt{1-R_v^2}$, then we
get $\Vert\ba_v\Vert= 1$.

Applying singular value decomposition (SVD) to $\bZ_B/\sqrt{n}$, we get
%
%
\begin{equation}\label{PF1_0}
\frac{1}{\sqrt{n}}\bZ_B = \bV\bM^{1/2} \bH^T,
\end{equation}
where $\bM= \operatorname{diag}(\mu_1,\ldots,\mu_{p-m})$ is a $(p-m) \times(p-m)$
diagonal matrix of ordered eigenvalues of
$\bS_{BB}$, $\bV$ is a $(p-m) \times(p-m)$ orthogonal matrix and
$\bH
$ is an $n \times(p-m)$ matrix. For $ n
\geq p-m$, $\bH$ has full rank orthogonal columns. When $ n < p-m $,
$\bH$ has more columns than rows, hence it
does not have full rank orthogonal columns. For the later case, we make
$\bH= [\bH_n,0]$ where $\bH_n$ is an $n
\times n$ orthogonal matrix.

\subsection{Propositions}

We introduce two propositions for later use. The proofs of the two
propositions can be found in Sections~\ref{L2_2} and~\ref{L2_3}.
\begin{Proposition}\label{prop1} Suppose $\bY$ is an $n \times m$
matrix with
fixed $m$ and each entry of
$\bY$ is i.i.d. random variable which satisfies the moment condition
of $z_{ij}$ in Assumption~\ref{assum1}. Let $ \bC$ be
an $ n \times n $ symmetric nonnegative definite random matrix and
independent of $\bY$. Further, assume $
\Vert\bC\Vert= O(1)$. Then
\[
\frac{1}{n} \bY^T \bC\bY- \frac{1}{n} \operatorname{trace}(\bC)\bI\stackrel
{p}{\rightarrow}0
\]
as $n \rightarrow
\infty$.\vadjust{\goodbreak}
\end{Proposition}
\begin{Proposition}\label{prop2}
Suppose $\mathbf{y}$ is an $ n $-dimensional random
vector which follows the same distribution of the row vectors of $\bY$
and independent of $\bS_{BB}$. Let $f(x)$ be a bounded continuous
function on
$[(1-\sqrt{\gamma})^2,(1+\sqrt{\gamma})^2]$ and $f(0) = 0$. Suppose
$\bF
=\mathrm{diag} (f(\mu_1),\ldots,f(\mu_{p-m}))$, where $\{ \mu_i
\}_{i=1}^{p-m}$ are ordered eigenvalues of $\bM$ which is defined on
(\ref{PF1_0}), then
\[
\frac{1}{n}\mathbf{y}^T \bH\bF\bH^T \mathbf{y}- \gamma\int{
f(x) \,dF_{\gamma} (x) }
\stackrel{p}{\rightarrow}0
\]
as $n \rightarrow\infty$, where $F_{\gamma} (x)$ is a distribution
function of Mar{\v{c}}enko--Pastur law with
parameter $\gamma$~\cite{marvcenko1967}.
\end{Proposition}

\subsection{\texorpdfstring{Proof of part \textup{(i)} of Lemma \protect\ref{lemma1}}
{Proof of part (i) of Lemma 1}}\label{L1}

\subsubsection{When $p$ is fixed} By the strong law of large numbers,
$\bS\stackrel{\mathrm{a.s.}}{\rightarrow}\bLambda$.
Since eigenvalues are continuous with respect to the
operator norm, the lemma follows after applying continuous mapping
theorem.

\subsubsection{\texorpdfstring{When $p \rightarrow\infty$}{When p converges to infinity}}

For every small $\varepsilon> 0$, there exist $\tilde{p}(n)$ and
$\gamma_{\varepsilon}$ such that $\tilde{p}(n) / n
\rightarrow\gamma_{\varepsilon} > 0$, $ \lambda_v(1+\gamma
_{\varepsilon
}/(\lambda_v - 1)) < \lambda_v + \varepsilon$ for all $v \leq m$,
$(1+\sqrt{\gamma_{\varepsilon}})^2 < 1 + \varepsilon$ and
$(1-\sqrt
{\gamma_{\varepsilon}})^2 > 1 - \varepsilon$. For
simplicity, we denote $\tilde{p}(n)$ as $\tilde{p}$. Suppose $\bZ
_{\tilde{p}}$ is a $ \tilde{p} \times n$ matrix
that satisfies the moment condition of $z_{ij}$ in Assumption \ref
{assum1}. Define
an augmented data matrix $\tilde{\bX}^T = [\bZ^T \bLambda,
\bZ_{\tilde{p}}^T]^T$ and its sample covariance matrix $\tilde{\bS} =
\tilde{\bX} \tilde{\bX}^T$. Let $\bS$ be a
$p\times p$ upper left submatrix of $\tilde{\bS}$. We also let $\hat
{\bS
}$ be an $(m+1)\times(m+1)$ upper left
submatrix of $\tilde{\bS}$. For $v \leq(m+1)$, by the interlacing
inequality (Theorem 4.3.15 of
Horn and Johnson~\cite{horn1990matrix}),
\[
\lambda_v(\hat{\bS}) \leq\lambda_v(\bS) \leq\lambda_v(\tilde
{\bS}).
\]
Since $\lambda_v(\hat{\bS}) \stackrel{\mathrm{a.s.}}{\rightarrow}\lambda
_v$, $\lambda_v(\tilde
{\bS})
\stackrel{\mathrm{a.s.}}{\rightarrow}
\lambda_v(1+\gamma_{\varepsilon}/(\lambda_v - 1)) < 1 + \varepsilon$
for $v \leq m$, and $\lambda_v(\tilde{\bS})
\stackrel{\mathrm{a.s.}}{\rightarrow}(1+\sqrt{\gamma_{\varepsilon}})^2 < 1 +
\varepsilon$ for $v =
m+1$, we have
\[
\lambda_v - o(1) \leq\lambda_v(\bS) < \lambda_v + \varepsilon+ o(1)
\qquad\mbox{for } v \leq m+1.
\]
Thus,
%
%
\begin{equation}\label{L1_1}
\lambda_v(\bS) \stackrel{\mathrm{a.s.}}{\rightarrow}\lambda_v\qquad
\mbox{for } v \leq m+1.
\end{equation}
Similarly by the interlacing inequality, we get
\[
\lambda_{\tilde{p}}(\tilde{S})\leq\lambda_p(S) \leq\lambda_{m+1}(S).
\]
Since $\lambda_{m+1}(S) \stackrel{\mathrm{a.s.}}{\rightarrow}1$ and $ \lambda
_{\tilde{p}}(\tilde{S})
\stackrel{\mathrm{a.s.}}{\rightarrow}(1-\sqrt{\gamma_{\varepsilon}})^2 >
1 - \varepsilon$, we conclude that
%
%
\begin{equation}\label{L1_2}
\lambda_p(S) \stackrel{\mathrm{a.s.}}{\rightarrow}1.
\end{equation}
 Part (i) of Lemma~\ref{lemma1} follows by (\ref{L1_1}) and
(\ref{L1_2}).

\subsection{\texorpdfstring{Proof of part \textup{(ii)} of Lemma \protect\ref{lemma2}}
{Proof of part (ii) of Lemma 2}}\label{L2}
Our proof of Lemma~\ref{lemma2}(ii) closely follows the arguments in
Paul \cite
{paul2007}.
From~\cite{paul2007}, it can be shown that
%
%
\begin{equation}\label{PF1_1}
\biggl(\bS_{AA} + \frac{1}{n}\bLambda_A^{1/2} \bZ_A \bH\bM(d_v \bI-
\bM
)^{-1} \bH^T \bZ_A^T \bLambda_A^{1/2} \biggr)\ba_v = d_v \ba_v
\end{equation}
and
%
%
\begin{equation}\label{PF1_2}
\ba_v^T\biggl(\bI+ \frac{1}{n}\bLambda_A^{1/2} \bZ_A \bH\bM(d_v \bI-
\bM
)^{-2} \bH^T \bZ_A^T \bLambda_A^{1/2} \biggr)\ba_v = \frac{1}{1-R_v^2},
\end{equation}
where $\bLambda_A = \mathrm{diag} \{ \lambda_1,\ldots,\lambda_m \} $.

\subsubsection{\texorpdfstring{When $\lambda_v > 1+\sqrt{\gamma}$}{When lambda v > 1 + sqrt gamma}}

We can show that
%
%
\begin{equation}\label{PF1_3}
\langle\ba_v,\mathbf{e}_{A,v} \rangle\stackrel{p}{\rightarrow}1
\end{equation}
and
%
%
\begin{equation}\label{PF1_4}
\frac{1}{n}\mathbf{z}_{Av}^T \bH\bM(d_v \bI- \bM)^{-2} \bH^T
\mathbf{z}_{Av}
\stackrel{p}{\rightarrow}
\cases{\displaystyle\gamma\int{\frac{x}{(\rho_v - x)^2}\,dF_{\gamma}(x)}, &\quad for
$\gamma> 0$,
\vspace*{2pt}\cr
0, &\quad for $\gamma= 0$,}\hspace*{-35pt}
\end{equation}
where $\mathbf{e}_{A,v}$ is a vector of the first $m$ coordinates of
the $v$th
population eigenvector $\mathbf{e}_v$,
$\rho_v$ is $ \lambda_v ( 1+
\frac{\gamma}{\lambda_v -1} )$  and $\mathbf{z}_{Av}$ is a\vspace*{2pt} vector
of $v$th row
of $\bZ_A$. The proofs can be
found in Section~\ref{L2_1}. Note that $\mathbf{e}_v$ is a vector with $1$ in its
$v$th coordinate and $0$ elsewhere. WLOG, we
assume that $ \langle\mathbf{e}_v,\mathbf{u}_v \rangle\geq0$.
Since $ \langle\mathbf{e}
_v,\mathbf{u}_v \rangle= \sqrt{1-R_v^2} \langle\mathbf
{e}_{A,v},\ba_v \rangle$,
$ \langle\mathbf{e}_v,\mathbf{u}_v \rangle\stackrel{p}{\rightarrow
}\sqrt{1-R_v^2}$. By
(\ref{PF1_2}) and (\ref{PF1_4}), we can show that
%
%
\begin{equation}\label{PF1_5}
\frac{1}{1-R_v^2} \stackrel{p}{\rightarrow}
\cases{\displaystyle 1 + \lambda_v \gamma\int{\frac{x}{(\rho_v - x)^2}\, d
F_{\gamma
}(x)}, &\quad for $\gamma> 0$,
\vspace*{2pt}\cr
1, &\quad for $\gamma= 0$.}
\end{equation}
From Lemma B.2 of~\cite{paul2007},
%
%
\begin{equation}\label{PF1_6}
\int{\frac{x}{(\rho_v - x)^2} \,d F_{\gamma}(x) }
= \frac{1}{(\lambda_v - 1)^2 - \gamma}.
\end{equation}
Thus,
%
%
\begin{equation}\label{P4}
\sqrt{1-R_v^2} \stackrel{p}{\rightarrow}
\cases{\displaystyle\sqrt{\biggl(1-\frac{\gamma}{(\lambda_v - 1)^2}\biggr) \Big/ \biggl(1 + \frac
{\gamma
}{\lambda_v - 1}\biggr)}, &\quad for $\gamma> 0$,
\vspace*{2pt}\cr
1, &\quad for $\gamma= 0$.}
\end{equation}
It concludes the proof of the first part of Lemma~\ref{lemma2}(ii).

\subsubsection{\texorpdfstring{When $1 < \lambda_v \leq1+\sqrt{\gamma}$}{When 1 < lambda v <= 1 + sqrt gamma}}

Here, we only need to consider $\gamma> 0$ because no eigenvalue
satisfies this condition when $\gamma=0$. We
first show that $R_v \stackrel{p}{\rightarrow}1 $, which implies
$\mathbf{u}_{A,v} \stackrel{p}{\rightarrow}0$,
hence\vadjust{\goodbreak} $ \langle\mathbf{e}_v,\mathbf{u}_v \rangle
\stackrel{p}{\rightarrow}0 $. For any $\varepsilon> 0$ and $x \geq
0$, define
\[
(x)_{\varepsilon} = \cases{x, &\quad if $x > \varepsilon$, \cr
\varepsilon, &\quad if $x \leq\varepsilon$,}
\]
and
\[
\bG_{\varepsilon} = \mathrm{diag}\bigl(d_v/\bigl((d_v - \mu_1)^2\bigr)_{\varepsilon
},\ldots,d_v/\bigl((d_v - \mu_{p-m})^2\bigr)_{\varepsilon}\bigr),
\]
then by Propositions~\ref{prop1} and~\ref{prop2},
%
%
\begin{equation}\label{PF12_1}
\frac{1}{n}\mathbf{z}_{Av}^T \bH\bG_{\varepsilon} \bH^T \mathbf
{z}_{Av} \stackrel{p}{\rightarrow}
\gamma\int{\frac{x}{((\rho_v - x)^2)_{\varepsilon}}\,dF_{\gamma}(x)}.
\end{equation}
By monotone convergence theorem,
%
%
\begin{equation}\label{PF12_2}
\gamma\int{\frac{x}{((\rho_v - x)^2)_{\varepsilon}}\,dF_{\gamma}(x)}
\stackrel{\varepsilon\rightarrow0 }{\longrightarrow}
\gamma\int{\frac{x}{(\rho_v - x)^2}\,dF_{\gamma}(x)}.
\end{equation}
The right-hand side of (\ref{PF12_2}) is
%
%
\begin{equation}\label{PF12_3}
\int_{a}^{b}{\frac{\sqrt{(b-x)(x-a)}}{2 \pi(\rho_v - x)^2}\,dx},
\end{equation}
where $a = (1-\sqrt{\gamma})^2$ and $b=(1+\sqrt{\gamma})^2$. Since
(\ref{PF12_3}) equals $\infty$ for any $a \leq
\rho_v \leq b$, we conclude that
%
%
\begin{equation}\label{PF12_4}
\frac{1}{n}\mathbf{z}_{Av}^T \bH\bM(d_v \bI- \bM)^{-2} \bH^T
\mathbf{z}_{Av} \stackrel{p}{\rightarrow}
\infty.
\end{equation}
Therefore $R_v \stackrel{p}{\rightarrow}1 $, which proves the second
part of Lemma~\ref{lemma2}(ii).

\subsubsection{\texorpdfstring{Proof of (\protect\ref{PF1_3}) and (\protect\ref{PF1_4})}
{Proof of (14) and (15)}}\label{L2_1}
Define
\begin{eqnarray*}
\mathcal{R}_v &=& \sum_{k \ne v}^{m} \frac{\lambda_v}{\rho_v(\lambda
_k -
\lambda_v)} \mathbf{e}_{A,k}
\mathbf{e}_{A,k}^T,\\
\mathcal{D}_v &=& \bS_{AA} + \bS_{AB}(d_v \bI-
\bS_{BB})^{-1}\bS_{BA} - (\rho_v / \lambda_v)\bLambda_{A} , \\
\alpha_v &=& \Vert\mathcal{R}_v \mathcal{D}_v\Vert+ | d_v - \rho_v |
\Vert\mathcal{R}_v\Vert\quad\mbox{and}\quad\beta_v = \Vert
\mathcal{R}_v \mathcal{D}_v \mathbf{e}_{A,v}\Vert.
\end{eqnarray*}

With the exactly same argument of~\cite{paul2007}, it can be shown that
\[
\ba_v - \mathbf{e}_{A,v} = - \mathcal{R}_v \mathcal{D}_v \mathbf
{e}_{A,v} + \br_v,
\]
where $ \br_v = -(1- \langle\mathbf{e}_{A,v},\ba_v \rangle)\mathbf
{e}_{A,v} - \mathcal
{R}_v \mathcal{D}_v (\ba_v - \mathbf{e}_{A,v}) +
(d_v - \rho_v) \mathcal{R}_v(\ba_v - \mathbf{e}_{A,v})$. By Lemma 1
of~\cite{paul2005}, $r_v = o_p(1)$, if $\alpha_v =
o_p(1)$ and $\beta_v = o_p(1)$.

When $\gamma=0$, $\bS_{AA} - (\rho_v / \lambda_v)\bLambda_{A}
\stackrel{p}{\rightarrow}0$
and the remainder of $\mathcal{D}_v$ is
%
%
\begin{equation}\label{PF13_1}\quad
\bS_{AB}(d_v \bI- \bS_{BB})^{-1}\bS_{BA}
= \frac{1}{n}\bLambda_A^{1/2} \bZ_A \bH\bM(d_v \bI- \bM)^{-1}
\bH^T
\bZ_A^T \bLambda_A^{1/2}.
\end{equation}
Since $d_v \stackrel{\mathrm{a.s.}}{\rightarrow}\lambda_v$ and $\mu_1
\stackrel{\mathrm{a.s.}}{\rightarrow}1$,
\[
\Vert\bH\bM(d_v \bI- \bM)^{-1} \bH^T \Vert\stackrel
{\mathrm{a.s.}}{\rightarrow}1/(\lambda_v -1).\vadjust{\goodbreak}
\]
By Proposition~\ref{prop1},
%
%
\begin{equation}\label{PF13_2}
0 \leq\Vert(\ref{PF13_1}) \Vert\leq\lambda_1 \frac{p \mu
_1}{n(d_v - \mu
_1)} +o_p(1) = o_p(1) ,
\end{equation}
hence $ \mathcal{D}_v = o_p(1) $.

When $\gamma> 0 $, $\mathcal{D}_v$ can be written as
%
%
\begin{eqnarray}\hspace*{28pt}
\mathcal{D}_v &=& [\bS_{AA} - \bLambda_{A} ]\nonumber\\
&&{}+ \biggl[\bLambda_{A}^{1/2}\biggl( \frac{1}{n}\bZ_A \bH\bM(\rho_v \bI- \bM
)^{-1} \bH^T \bZ_A\nonumber\\
&&\hspace*{49pt}{} - \frac{1}{n}\operatorname{trace}\bigl(\bM(\rho_v \bI- \bM)^{-1}\bigr)
\bI
\biggr)\bLambda_{A}^{1/2} \biggr] \\
&&{}+ \biggl[\biggl(\frac{1}{n}\operatorname{trace}\bigl( \bM(\rho_v \bI- \bM)^{-1}\bigr) - \gamma\int
{\frac
{x}{\rho_v -x} \,d F_{\gamma} (x)}\biggr) \bLambda_{A}\biggr] \nonumber\\
&&{}+ \biggl[(\rho_v - d_v) \frac{1}{n}\bLambda_{A}^{1/2}\bZ_A \bH\bM(\rho_v
\bI- \bM)^{-1}(d_v \bI- \bM)^{-1} \bH^T \bZ_A
\bLambda_{A}^{1/2}\biggr].\nonumber
\end{eqnarray}
The first term of the right-hand side is $o_p(1)$ by the weak law of large number. The
second and third terms are $o_p(1)$ by
Propositions~\ref{prop1} and~\ref{prop2}. For the fourth term, $\rho_v
- d_v = o_p(1)$ and
its remainder part is $O_p(1)$.
Therefore, $\mathcal{D}_v = o_p(1)$.
By combining the above results and
$\mathcal{R}_v = O_p(1)$ plus $d_v - \rho_v = o_p(1)$, we prove~(\ref{PF1_3}).

For (\ref{PF1_4}): When $\gamma= 0 $, (\ref{PF1_4}) can be proved by
the exactly same way used
to show (\ref{PF13_2}). When $\gamma> 0$, $d_v \stackrel
{\mathrm{a.s.}}{\rightarrow}\rho_v$, and $
\mu_1 \stackrel{\mathrm{a.s.}}{\rightarrow}(1+ \sqrt{\gamma})^2 < \rho_v
$, hence $ \Vert\bC\Vert\stackrel{\mathrm{a.s.}}{\rightarrow}\frac{(1+
\sqrt{\gamma})^2 }{(\rho
_v - (1+
\sqrt{\gamma})^2 )^2} $. Therefore, the
result follows according to Propositions~\ref{prop1} and~\ref{prop2}.

\subsection{\texorpdfstring{Proof of Proposition \protect\ref{prop1}}{Proof of Proposition 1}}\label{L2_2}
Let $\mu_1 \geq\mu_2 \geq\cdots\geq\mu_n $ be the ordered
eigenvalues of $\mathbf{C}$, and $c_{ij}$ be the
$(i,j)$th element of $\bC$. Suppose $\mathbf{y}_s$ is the $s$th
column of $\bY
$, and $y_{ij}$ is the $(i,j)$th element
of $\bY$. We further define $ \psi(s,s) = \frac{1}{n} \mathbf
{y}_s^T \bC\mathbf{y}_s
- \frac{1}{n} \operatorname{trace} (\bC) $ and $
\psi(s,t) = \frac{1}{n} \mathbf{y}_s^T \bC\mathbf{y}_t $ for $s
\ne t$. The
conditional mean of $\psi(s,s)$ given $\bC$ is
%
%
\begin{eqnarray}
E( \psi(s,s) |\bC)
&=& E\biggl( \frac{1}{n} \sum_{i,j} c_{ij} y_{is} y_{js} |\bC\biggr) - \frac{1}{n}
\sum_{i=1}^n \mu_i \nonumber\\
&=& \frac{1}{n} \sum_{i=1}^n c_{ii}E( y_{is}^2 ) + \frac{2}{n} \sum
_{i<j}^n c_{ij}E( y_{is} y_{js} ) - \frac{1}{n} \sum_{i=1}^n \mu_i
\\
&=& \frac{1}{n} \sum_{i=1}^n c_{ii} - \frac{1}{n} \sum_{i=1}^n \mu
_i =0.\nonumber
\end{eqnarray}
Thus, $ E( \psi(s,s) ) = E(E( \psi(s,s) |\bC)) = E(0) = 0$.\vadjust{\goodbreak}

Next, the conditional variance of $\psi(s,s)$ given $\bC$ is
%
%
\begin{eqnarray}\label{P_1_2}
\operatorname{Var}(\psi(s,s) | \bC) &=& \frac{1}{n^2} \operatorname{Var}\biggl(\sum_{i,j} c_{ij} y_{is}
y_{js} | \bC\biggr)\nonumber\\
&=& \frac{1}{n^2} \sum_{i,j,l,q=1}^n c_{ij} c_{lq} \operatorname{Cov}(y_{is} y_{js},
y_{ls} y_{qs})\nonumber\\[-8pt]\\[-8pt]
&=& \frac{4}{n^2} \sum_{i,j=1}^n c_{ij}^2 \operatorname{Var}(y_{is} y_{js}) \nonumber
\leq\frac{4 \alpha}{n^2} \sum_{i,j=1}^n c_{ij}^2
\nonumber\\
&=& \frac{4 \alpha}{n^2} \operatorname{trace}( \bC^2)
= \frac{4 \alpha}{n^2} \sum_{i=1}^n \mu_i^2,\nonumber
\end{eqnarray}
where $\alpha= \max(1,E(y_{is}^4)-1)$.
Since $\Vert\bC\Vert= O(1)$, $\mu_i^2 \leq\Vert\bC\Vert^2 =
O(1)$. Therefore,
$\operatorname{Var}(\psi(s,s) | \bC) \leq O(1/n)$ and $ \operatorname{Var}(\psi(s,s)) = \operatorname{Var}(E(\psi
(s,s)|\bC))+\break E(\operatorname{Var}(\psi(s,s)|\bC)) \leq0 +
O(1/n) \rightarrow0$ as $n \rightarrow\infty$. By the Chebyshev
inequality, we can conclude that
\[
\psi(s,s) \stackrel{p}{\rightarrow} 0.
\]
We can similarly show $ \psi(s,t) \stackrel{p}{\rightarrow} 0 $, which
we omit here.

\subsection{\texorpdfstring{Proof of Proposition \protect\ref{prop2}}{Proof of Proposition 2}}\label{L2_3}
Consider an expansion
\begin{eqnarray*}
&& \frac{1}{n}y^T \mathbf{H} \mathbf{F} \mathbf{H}^T y - \gamma
\int{
f(x) \,dF_{\gamma} (x) } \\
&&\qquad = \biggl[\frac{1}{n}y^T \mathbf{H} \mathbf{F} \mathbf{H}^T y - \frac
{1}{n}\operatorname{trace}( \mathbf{F} ) \biggr] \\
&&\qquad\quad{} + \biggl[ \frac{1}{n}\operatorname{trace}( \mathbf{F} ) - \gamma\int{ f(x) \,dF_{\gamma}
(x) }\biggr] \\
&&\qquad = \mbox{(a)} + \mbox{(b)}.
\end{eqnarray*}

We show that both (a) and (b) converge to 0 in probability.

(a): Since $\mu_1 \stackrel{\mathrm{a.s.}}{\rightarrow} (1+\sqrt{\gamma})^2 $,
$\mu_{\min(p-m,n)} \stackrel{\mathrm{a.s.}}{\rightarrow} (1-\sqrt{\gamma})^2$,
$\mu_k =0$ for $k > \min(p-m,n)$  and $f(x)$ is continuous and bounded
on $[(1-\sqrt{\gamma})^2,(1+\sqrt{\gamma})^2]$, there
exists $ K > 0$ such that
$ {\sup_{i}} |f(\mu_i)| < K$ a.s.
Let $ \mathbf{C} = \mathbf{H}
\mathbf{F} \mathbf{H}^T $, then $ \operatorname{trace} (\mathbf{C}) = \operatorname{trace}(\mathbf
{F}) $. By Proposition~\ref{prop1}, $ \mbox{(a)} = o_p(1) $.

(b): Let $F_{p-m}$ be an empirical spectral distribution of $\bS
_{BB}$, then
\[
\frac{1}{n} \operatorname{trace}( \mathbf{F}) = \frac{p-m}{n}\int{ f(x) \,dF_{p-m}
(x) }
\]
and $ \int{ f(x) \,dF_{n} (x)} \stackrel{p}{\rightarrow} \int{ f(x)
\,dF_{\gamma} (x)} $
\cite{marvcenko1967,bai1999msa}. Thus,
\[
\frac{p-m}{n}\int{ f(x) \,dF_{p-m} (x) } \stackrel{p}{\rightarrow}
\gamma\int{ f(x) \,dF_{\gamma} (x) },
\]
which shows that $\mbox{(b)} = o_p(1)$.

Combining (a) and (b), we finish the proof.

\subsection{\texorpdfstring{Proof of Theorem \protect\ref{theo1}}{Proof of Theorem 1}}\label{T1}
Without loss of generality, we assume
$ \langle\mathbf{g}_v,\tilde{\mathbf{p}}_v \rangle\geq0$. Let
$\mathbf{e}_v =
\{\mathbf{e}_{A,v}, \mathbf{e}_{B,v} \} $, then $\mathbf{e}_{A,v}$ is
the vector with 1 in
$v$th coordinate and~$0$ elsewhere, and
$\mathbf{e}_{B,v}$ is the zero vector. Since $ \bS_{AA} \mathbf
{u}_{A,v} + \bS_{AB}
\mathbf{u}_{B,v} = d_v \mathbf{u}_{A,v}$, we have
%
%
\begin{eqnarray}\label{PF3_1}
\langle \mathbf{g}_v,\tilde{\mathbf{p}}_v\rangle &=& \frac{1}{n \sqrt{d_v
\lambda_v} } \mathbf{e}
_v^T \bX\bX^T \mathbf{u}_v \nonumber\\
&=& \mathbf{e}_{A,v}^T \bS_{AA} \mathbf{u}_{A,v} / \sqrt{d_v \lambda
_v} + \mathbf{e}_{A,v}^T
\bS_{AB} \mathbf{u}_{B,v} / \sqrt{d_v \lambda_v}
\nonumber\\[-8pt]\\[-8pt]
&=& \frac{d_v}{\sqrt{d_v \lambda_v} } \mathbf{e}_{A,v}^T \mathbf
{u}_{A,v} = \sqrt
{\frac{d_v}{\lambda_v}} \mathbf{e}_v^T \mathbf{u}_v \nonumber\\
&\stackrel{p}{\rightarrow}&
\cases{\sqrt{\biggl(1-\dfrac{\gamma}{(\lambda_v-1)^2}\biggr)}, &\quad for $\lambda
_v >
1+\sqrt{\gamma}$,
\vspace*{2pt}\cr
0, &\quad for $1 < \lambda_v \leq1+\sqrt{\gamma}$.}\nonumber
\end{eqnarray}

\subsection{\texorpdfstring{Proof of Theorem \protect\ref{theo2}}{Proof of Theorem 2}}\label{T2}

First, we show the square of the denominator converges to $\rho
(\lambda_v)$.
Since $ p_{vj} = \mathbf{u}_v^T \mathbf{x}_j$, and $E(p_{vi}^2) =
E(p_{vj}^2)$ for
$i\neq j$,
%
%
\begin{eqnarray}\label{PF2_0}
E(p_{vj}^2) & = &\frac{1}{n}E\Biggl(\sum_{j=1}^{n}p_{vj}^2\Biggr)
=\frac{1}{n}E\Biggl(\sum_{j=1}^{n} (\mathbf{u}_v^T \mathbf{x}_j)^2\Biggr)
\nonumber\\[-8pt]\\[-8pt]
& = &E(\mathbf{u}_v^T \bX\bX^T \mathbf{u}_v / n)
= E(d_v) \stackrel{\mathrm{a.s.}}{\rightarrow}\rho(\lambda_v).\nonumber
\end{eqnarray}
Next, we show the square of numerator converges to $\phi(\lambda_v)^2
(\lambda_v -1 ) + 1$.
Define $\mathbf{u}_v^{\bot} := \frac{1}{\sqrt{1-(\mathbf
{u}_v^T\mathbf{e}_v)^2}}(I-\mathbf{e}_v \mathbf{e}
_v^T)\mathbf{u}_v$, then $\mathbf{u}_v$ can be
expressed as
\[
\mathbf{u}_v = (\mathbf{u}_v^T\mathbf{e}_v) \mathbf{e}_v + \sqrt{1
- (\mathbf{u}_v^T\mathbf{e}_v)^2 } \mathbf{u}_v^{\bot}.
\]
Partition $ \mathbf{u}_v^{\bot} = \{ \mathbf{u}^{\bot}_{A,v} ,
\mathbf{u}^{\bot}_{B,v} \}$.
From (\ref{PF1_3}),
$\mathbf{a}_v \stackrel{p}{\rightarrow}\mathbf{e}_{A,v} $,
therefore $\mathbf{u}^{\bot
}_{A,v} \stackrel
{p}{\rightarrow} 0$ and $ \mathbf{u}^{\bot
T}_{B,v} \mathbf{u}^{\bot}_{B,v} \stackrel{p}{\rightarrow} 1 $.
Since $\mathbf{x}
_{\mathrm{new}} $ and $\mathbf{u}_v$ are independent, we have
%
%
\begin{eqnarray}\label{PF2_1}
E(q_v^2|\mathbf{u}_v) & = & E( (\mathbf{u}_v^T \mathbf{x}_{\mathrm{new}})^2 |
\mathbf{u}_v) = \mathbf{u}_v^T E( \mathbf{x}
_{\mathrm{new}} \mathbf{x}_{\mathrm{new}}^T | \mathbf{u}_v) \mathbf{u}_v
= \mathbf{u}_v^T \bolds\Lambda\mathbf{u}_v \nonumber\\
&=& (\mathbf{u}_v^T\mathbf{e}_v)^2 \mathbf{e}_v^T \bolds\Lambda\mathbf
{e}_v + \bigl(1 - (\mathbf{u}_v^T\mathbf{e}_v)^2 \bigr) \mathbf{u}
_v^{\bot T} \bolds\Lambda\mathbf{u}_v^{\bot}\nonumber\\
&&{} + 2 \mathbf{u}_v^T\mathbf{e}_v \sqrt{1 - (\mathbf{u}_v^T\mathbf
{e}_v)^2} \mathbf{e}_v^T \bolds\Lambda\mathbf{u}_v^{\bot
} \nonumber\\[-8pt]\\[-8pt]
&=& (\mathbf{u}_v^T\mathbf{e}_v)^2 \lambda_v + \bigl(1 - (\mathbf
{u}_v^T\mathbf{e}_v)^2 \bigr) ( \mathbf{u}^{\bot T}_{A,v}
\bolds\Lambda_A \mathbf{u}^{\bot}_{A,v} + \mathbf{u}^{\bot T}_{B,v}
\mathbf{u}^{\bot}_{B,v} )
\nonumber\\
&&{}+ 2 \mathbf{u}_v^T\mathbf{e}_v \sqrt{1 - (\mathbf{u}_v^T\mathbf
{e}_v)^2} \mathbf{e}_{A,v} \bolds\Lambda_A \mathbf{u}
^{\bot}_{A,v} \nonumber\\
&\stackrel{p}{\rightarrow}&\phi(\lambda_v)^2 (\lambda_v -1 ) +
1.\nonumber
\end{eqnarray}
From (\ref{PF2_0}) and (\ref{PF2_1}),
%
%
\begin{equation}\label{PF2_3}
\sqrt{\frac{E(q_v^2)}{E(p_{vi}^2)}} \rightarrow
\sqrt{\frac{\phi(\lambda_v)^2 (\lambda_v -1 ) + 1}{\rho(\lambda_v)}}
=\frac{(\lambda_v -1)}{(\lambda_v +\gamma-1)}.
\end{equation}

\subsection{\texorpdfstring{Proof of Theorem \protect\ref{theo3}}{Proof of Theorem 3}}\label{T3}

Since $ \rho^{-1}(p r_v) \rightarrow\lambda_v $ for $ v \leq k$, WLOG
we assume that $k_0 = k$, where $k$ is
the number of $\lambda_v$ bigger than $ 1+ \sqrt{\gamma}$. Set
%
%
\begin{equation}\label{PF4_1}
h(x) = \sum_{v=1}^{k} \rho^{-1}( r_v x) + p - k -x.
\end{equation}
The first and second partial derivatives of $h(x)$ are
%
%
\begin{eqnarray}\label{PF4_2}
\frac{ \partial h(x)}{\partial x} &=& \frac{1}{2} \sum_{v=1}^{k} r_v
+ \frac{1}{2} \sum_{v=1}^{k} \frac{ (x r_v -(1+\gamma)) r_v}{\sqrt{(x
r_v -(1+\gamma))^2 - 4\gamma}}
-1,
\\
%
%
\label{PF4_3}
\frac{ \partial^2 h(x)}{\partial x^2} &=&
2\sum_{v=1}^{k} \frac{-r_v^2 \gamma}{((x r_v -(1+\gamma))^2 -
4\gamma
)^{3/2} }<0,
\end{eqnarray}
so $h(x)$ is a concave function of $x$ given $r_v$. From the fact that
$\rho^{-1}(r_v p) > 1$ for $v \leq k$, we
know $h(p) > 0 $. Because of the concave nature of this function,
$h(x)=0$ has a unique solution $\tau$ on $
[p, \infty)$, which $\sum_{v=1}^{k_l} \hat{\lambda}_{v,l} + p - m_l$
converges to. Thus, $ \hat{d}_v =
\tau r_v. $ Define $ \tilde{d}_v = r_v \omega$ where $ \omega= \sum
_{v=1}^{k} \lambda_v + p -k $, and
set $d_v$ as the sample eigenvalue when $\sigma^2 = 1$. The sum of all
$d_v$ is
%
%
\begin{equation}\label{PF4_6}
\sum_{v=1}^{p} d_v = \frac{1}{n} \operatorname{trace}(\mathbf{Z} \mathbf{Z}^T \bolds\Lambda)
= \frac{1}{n} \sum_{i=1}^{p} \sum_{j=1}^{n} \lambda_i z_{ij}^2,
\end{equation}
thus
%
%
\begin{equation}\label{PF4_7}
E \biggl( \frac{\sum_{v=1}^{p} d_v}{\omega} \biggr) = \frac{\sum_{v=1}^{m}
\lambda
_v + p -k}{\omega}
\rightarrow1
\end{equation}
and
%
%
\begin{equation}\label{PF4_8}
\operatorname{Var} \biggl( \frac{\sum_{v=1}^{p} d_v}{\omega} \biggr) =
\frac{1}{n}\frac{\sum_{v=1}^{p} \lambda_v^2}{\omega^2}
\bigl(E(z_{11}^4) -1\bigr) \rightarrow0.
\end{equation}
By (\ref{PF4_7}) and (\ref{PF4_8}),
%
%
\begin{equation}\label{PF4_9}
\sum_{v=1}^{p} d_v / \omega= 1 + o_p(1).\vadjust{\goodbreak}
\end{equation}
Since $d_v \rightarrow\rho(\lambda_v) $ for $ v \leq k$,
%
%
\begin{equation}\label{PF4_10}
\tilde{d}_v = d_v \omega\Big/ \sum_{v=1}^{p} d_v = d_v\bigl( 1+ o_p(1)\bigr)
\stackrel{p}{\rightarrow} \rho(\lambda_v).
\end{equation}
Now, we show that $ \tau= \omega+ o_p(1) $. Plugging $\omega$ into
$h(x)$ and combining the fact that
$\rho^{-1}( \tilde{d}_v) = \lambda_v + o_p(1)$, we get
%
%
\begin{equation}\label{PF4_11}
h(\omega)
= \sum_{v=1}^{k} \rho^{-1}( \tilde{d}_v) - \sum_{v=1}^{k} \lambda
_v = o_p(1).
\end{equation}
From the facts that $h(x)$ is a continuous concave function, $ \omega
>p $, and $h(p) > 0 $, we conclude that
%
%
\begin{equation}\label{PF4_13}
\omega= \tau+ o_p(1).
\end{equation}
Therefore,
%
%
\begin{equation}\label{PF4_14}
\hat{d}_v = r_v \tau= r_v \bigl(\omega+ o_p(1)\bigr) = \tilde{d}_v + o_p(1)
\stackrel{p}{\rightarrow} \rho(\lambda_v)
\end{equation}
for $ v \leq k $, which concludes the proof.

\printaddresses

\end{document}